\documentstyle{article}

\begin{document}

\title{
Spaces of holomorphic almost periodic functions
on a strip}
\author{Favorov S.Yu., Udodova O.I.}
\date{}
\maketitle

\centerline
{\begin{tabular}{c}
Karazin Kharkov National University\\
Department of Mechanic and Mathematic\\
Sq. Svobody 4, Kharkov, 61077, Ukraine\\
E-mail: favorov@assa.vl.net.ua,\
Olga.I.Udodova@univer.kharkov.ua\\
\end{tabular}}

\bigskip{
The notions of almost periodicity in the sense of Weyl and Besicovitch
of the order $p \ge 1$ are extended  to
holomorphic functions on a strip.
We prove that the spaces of holomorphic almost periodic functions
in the sense of Weyl for various orders
 $p$ are the same. These spaces are considerably wider than the space of holomorphic
uniformly almost periodic functions and considerably narrower than the spaces
of holomorphic almost periodic functions in the sense of Besicovitch.
Besides we construct examples showing that the spaces of holomorphic
almost periodic functions in the sense of Besicovitch  for various orders
 $p$  are all different.\\
{ \it 2000 Mathematics Subject Classification 42A75, 30B50.}}

\bigskip

A continuous function $f$ on a strip $\Pi _{[ {a,b}]}=
\left\{ {z = x + iy \in {\bf C} : a \le y \le b} \right\}$ is
called almost periodic if
for any
$\varepsilon  > 0$ the set of $\varepsilon $-almost periods
$$\left\{ {\tau  \in {\bf R}:d_{[a,b]}^U (f(z),f(z + \tau ))
 < \varepsilon } \right\}$$
is relatively dense on ${\bf R}$, i.e.,
 its intersection  with any segment of the length
$L=L(\varepsilon)$ is nonempty.
Here
\begin{equation}
\label{d_U}
d_{[a,b]}^U (g,h) = \mathop {\sup }\limits_{z \in \Pi _{\left[ {a,b}
 \right]} } |g(z) -
h(z)|
\end{equation}
is the standard uniform metric on $\Pi_{[a,b]}$.

We denote by $U_{[a,b]} $ the space of almost periodic functions
on $\Pi_{[a,b]}$.
By the Approximation Theorem, $U_{[a,b]}$ is equal to
the closure of the set of all finite exponential sums
\begin{equation}
\label{f1}
\sum\limits_{n = 1}^N {c_n (y)e^{i\lambda _n x} },\quad
 \lambda _n  \in {\bf R}, \quad  c_n(y)\in C_{[a,b]}
\end{equation}
 with respect to the metric $d_{[a,b]}^U $.

Note that any continuous periodic function on a strip
$\Pi _{\left[ {a,b} \right]} $ with real
periods is almost periodic
since any period is an $\varepsilon $-almost period for every $\varepsilon
> 0$, and the set of periods forms a dual-sided arithmetical progression.

We can replace the uniform metric $d_{[a,b]}^U $
by either the Stepanov metric
\begin{equation}
\label{d_S}
d_{[a,b]}^{S^p}(g,h) =
\mathop
{\sup}\limits_{z \in\Pi_{\left[ a,b \right]}}\left(
{\int\limits_0^1 {|g(z + t) - h(z + t)|^p dt}} \right)^\frac 1p,\;p \ge 1
\end{equation}
or the Weyl metric
\begin{equation}
\label{d_W}
d_{[a,b]}^{W^p } (g,h) = \mathop {\overline {\lim } }\limits_{T \to \infty } \mathop
{\sup }\limits_{z \in \Pi _{\left[ {a,b} \right]} } \left( {\frac{1}{{2T}}\int\limits_{ - T}^T {|g(z
+ t) - h(z + t)|^p dt} } \right)^{\frac{1}{p}} ,\;p \ge 1
\end{equation}
or the Besicovitch metric
\begin{equation}
\label{d_B}
d_{[a,b]}^{B^p } (g,h) = \mathop {\overline {\lim } }\limits_{T \to \infty } \mathop
{\sup }\limits_{a \le y \le b} \left( {\frac{1}{{2T}}\int\limits_{ - T}^T {|g(t + iy) - h(t + iy)|^p
dt} } \right)^{\frac{1}{p}} ,\;p \ge 1.
\end{equation}
We suppose that the functions $g(z)$ and $h(z)$ are measurable
and $|g(x+iy)|^p$ and $|h(x+iy)|^p$ are locally integrable with respect to
the variable $x$ for fixed $y$.

The closure of the set of sums (\ref{f1})
in metrics (\ref{d_S}) - (\ref{d_B})will be called
the space of almost periodic functions in the sense of
Stepanov of order  $p$,
the space of almost periodic functions
in the sense of Weyl of order $p$,
the space of almost periodic functions
in the sense of Besicovitch of order $p$,
and , respectively, will be denote by
\begin{equation}
\label{f2}
S_{[a,b]}^{^p } , W_{[a,b]}^{^p } \ {\rm and} \ B_{[a,b]}^{^p }.
\end{equation}

In the case $a = b = 0$, i.e., for functions defined on the real
axis, these spaces are well-known (see [2], [3], [6], p. 189-247.)
Note that we can also define the spaces (\ref{f2}) by using the concept of
$\varepsilon$-almost period. But for the space of
$B_{[a,b]}$
this definition is much more complicated
(see. [2], p. 91-104).

It is clear that every sum (\ref{f1}) is uniformly bounded on $\Pi _{\left[ {a,b}
\right]} $, and therefore an almost periodic function
$f(z)$ from the space $U_{[a,b]} $ (or from spaces (\ref{f2}))
is bounded in the corresponding metric, i.e.
$d_{[a,b]}^U (f,0) < \infty $ (or
$d_{[a,b]}^{S^p } (f,0) < \infty $,
$d_{[a,b]}^{W^p } (f,0) < \infty $, $d_{[a,b]}^{B^p } (f,0) < \infty $).
Observe that in the definitions of spaces
(\ref{f2}) we can replace  sums (\ref{f1}) by
functions from $U_{[a,b]} $.

It is easy to see that the elements of metric spaces (\ref{f2})
are equivalent classes of functions with the corresponding zero distances.
For example the equivalent class for the space $S_{[a,b]}^p $ consists of the functions
which coincide  a. e. on every horizontal straight line.
Note that every two functions with the difference
 $e^{ - \left( {x + iy} \right)^2 } $ always belong to the same class
in the spaces $W_{[a,b]}^p $ and $B_{[a,b]}^p $,
since $d_{[a,b]}^{W^p } (e^{ - \left(
{x + iy} \right)^2 } ,0) = 0$ for all $p < \infty $ and $a \le b$.
Nevertheless we will use the notation  $f \in W_{[a,b]}^p $
(or $f \in B_{[a,b]}^p $); this means that the equivalent class
of the function $f$ belongs to the corresponding space.

The H\"older inequality implies that for all $p < p'$
$$d_{[a,b]}^{S^p } (g,h) \le d_{[a,b]}^{S^{p'} } (g,h), \ d_{[a,b]}^{W^p } (g,h) \le
d_{[a,b]}^{W^{p'} } (g,h), \ d_{[a,b]}^{B^p } (g,h) \le d_{[a,b]}^{B^{p'} } (g,h).$$
Besides for all $p \ge 1$
$$d_{[a,b]}^{B^p } (g,h) \le d_{[a,b]}^{W^p } (g,h) \le d_{[a,b]}^{S^p } (g,h) \le
d_{[a,b]}^U (g,h)$$
(the inequality between $d_{[a,b]}^{W^p } (g,h)$ and $d_{[a,b]}^{S^p } (g,h)$
 follows from the relation
$$\frac{1}{{2L}}\int\limits_{ - L}^L {|f(x + t) - g(x + t)|^p dt}  \le
\frac{1}{{2L}}\sum\limits_{j =  - \left[ L \right] - 1}^{\left[ L \right]} {\int\limits_j^{j + 1}
{|f(x + t) - g(x + t)|^p dt}  \le \frac{{\left[ L \right] + 1}}{L}d_{\left[ {a,b} \right]}^{S^p } }
(f,g))$$
Hence we have
\begin{equation}
\label{f3}
S_{[a,b]}^{p'}  \subset S_{[a,b]}^p
, W_{[a,b]}^{p'}  \subset W_{[a,b]}^p
, B_{[a,b]}^{p'}  \subset B_{[a,b]}^p
\end{equation}
and
\begin{equation}
\label{f4}
U_{[a,b]}^{}  \subset S_{[a,b]}^p  \subset W_{[a,b]}^p  \subset B_{[a,b]}^p.
\end{equation}
Here all the inclusions mean that each equivalent class in a "narrower"
space is contained in some equivalent class in a "wider" space.
All the inclusions are strict even in the case $a = b = 0$:
there exists the equivalent class  in a "wider" space
 containing no class from a "narrower" space (see [4]).

Note that  there exists the mean value
\begin{equation}
\label{M_f}
\left( {M_t f} \right)(y) =
 \mathop {\lim}\limits_{T \to \infty } \frac{1}{{2T}}\int\limits_{ - T}^T
\,f(t,y)dy
\end{equation}
uniformly in $y \in \left[ {a,b} \right]$
for the before mentioned functions. Indeed,
the mean value of finite exponential sums (\ref{f1})
equals the coefficient at $e^{i0t} $;
in the general case, the existence of limits (\ref{M_f}) follows easily from
the definition of almost periodicity.
By the same way, we can prove  the relation
\begin{equation}
\label{f5}
\mathop {\lim }\limits_{T \to \infty } \left|
{\mathop {\sup }\limits_{\scriptstyle x \in {\bf R}^m ,
\hfill \atop
  \scriptstyle a \le y \le b \hfill}
\frac{1}{{2T}}\int\limits_{ - T}^T {f(x + t,y)dt - (M_t f)(y)}
} \right| = 0
\end{equation}
for the space $W^1 $ and for all "narrower" spaces.

Let $\Pi _{(a,b)}  = \left\{ {z = x + iy \in {\bf C}:-\infty\le a < y < b\le\infty} \right\}$
be an open strip, maybe infinite width.
A function $f(z)$ is uniformly almost periodic on
$\Pi _{(a,b)} $
(or almost periodic
in the sense of Stepanov, Weyl, and Besicovitch),
if the restriction $f(z)$ to every strip $\Pi _{\left[ {\alpha ,\beta }
\right]} $ for $a < \alpha
< \beta  < b$ is almost periodic in metrics (\ref{d_U}), (\ref{d_S}) - (\ref{d_B}).
We denote by
$HU_{(a,b)} $, $HS_{\left(
{a,b} \right)}^p $, $HW_{\left(
{a,b} \right)}^p $, $HB_{\left( {a,b} \right)}^p $
the corresponding spaces of holomorphic almost
periodic functions on $\Pi_{(a,b)}$ .
The inclusions similar to (\ref{f3}) and (\ref{f4})
hold for these spaces.
The spaces $HB_{(a,b)}^p $ were studied earlier in [7].
These spaces were also defined in [1]
as sets of holomorphic functions on a strip with the following property:
 the restrictions
to each straight line in the strip are almost periodic functions
in the sense of Besicovitch  of order $p$.
However as it follows from [1], theorem 3.4,
for functions growing as $O(e^{e^{|z|^N } } )$ on a strip
these definitions coincide.

The following Linfoot's theorem is well-known:

{\bf Theorem L.}
 (see. [2], p. 146)

{\it Spaces $HU_{\left( a,b \right)}$
and $HS_{\left( {a,b} \right)}^p$ coincide for all $p \ge 1$.}

Here we obtain the similar result:

{\bf Theorem 1.}

{\it The spaces $HW_{(a,b)}^p $ coincide for all $p \ge 1$.}

The proof of this theorem is based on the following proposition.

{\bf Proposition 1.}

{\it Each function $f \in HW_{(a,b)}^1 $ is uniformly bounded
on every substrip
$\Pi_{\left[ \alpha ,\beta \right]}$, $a < \alpha < \beta  < b$.}

Next, we prove that inclusions
 $HU_{\left( {a,b} \right)}^{}  \subset HW_{\left(
{a,b} \right)}^1 $, $HW_{\left( {a,b} \right)}^1
\subset HB_{\left( {a,b} \right)}^p $,
$HB_{\left( {a,b} \right)}^{p'}  \subset HB_{\left( {a,b} \right)}^p $,
$p' > p \ge 1$, are strict in the same sense as inclusions (\ref{f3}) and
(\ref{f4})).

{\bf Theorem 2.}

{\it There exists a function
$f \in HW_{\left( { - \infty ,\infty } \right)}^1 $ such that
every function $g$ equivalent to  $f$  in any space
$HW_{\left[ - H,H \right]}^1 $, does not belong to
 $HU_{\left[ - H,H \right]}$ for all $H>0$.}

{\bf Theorem 3.}

{\it There exists a function
 $f \in \bigcap\limits_{p \ge 1}^{} {HB_{\left( { - \infty ,\infty }
\right)}^p } $ such that
every function $g$ equivalent to  $f$  in any space
$HB_{\left[ - H,H \right]}^p $, $p \ge 1$, does not belong to
$HW_{\left[ - H,H \right]}^1 $ for all $H>0$.}

{\bf Theorem 4.}

{\it For all $p' > p \ge 1$ there exists a function
$f \in HB_{\left( { - \infty ,\infty } \right)}^p $
such that
every function $g$ equivalent to  $f$  in any space
$HB_{\left[ - H,H \right]}^p $,
does not belong to $HB_{\left[ - H,H \right]}^{p'} $
for all $H>0$.}

{\bf The proof of Proposition 1.}

Suppose $a < \alpha ' < \alpha  < \beta  < \beta ' < b$.
 Since $d_{\left[ {\alpha ',\beta '}
\right]}^{W^1 } (f,0) < \infty $,
we have
$$\frac{1}{{2T_0 }}\int\limits_{ - T_0 }^{T_0 } {|f(z + u)|du}  \le C$$
for some
 $C < \infty ,\;T_0  < \infty $ and all $z \in
\Pi _{[\alpha ',\beta ']} $.

Fix
 $r < \mathop {\min }\limits_{} \left\{ {T_0 ,\;\alpha  -
 \alpha ',\;\beta ' - \beta }
\right\}$.

Since the function $f(z)$ is holomorphic, we have
for all $z_0  \in \Pi _{[\alpha ,\beta ]} $
$$|f(z_0 )| \le \frac{1}{{\pi r^2 }}\int\limits_{|z - z_0 | < r}^{} {|f(z)|dxdy} $$
and
$$
 |f(z_0 )| \le \frac{1}{{\pi r^2 }}\int\limits_{ - T_0 }^{T_0 } {\int\limits_{|v| < r}^{}
{|f(z_0  + u + iv)|dudv \le } }  
 \frac{1}{{\pi r^2 }}\mathop {\sup }\limits_{z \in \Pi _{[\alpha ,\beta ]} } \int\limits_{ -
T_0 }^{T_0 } {|f(z + u)|du}  \le \frac{{2T_0 C}}{{\pi r^2 }}. \\
$$
Proposition 1 is proved.

Consider the sums
$$K^{(q)} (t) = \sum\limits_n^{} {k_n^q e^{ - i\lambda _n t} },$$
which are called Bochner-Fejer kernels
(see,
for example, [6], p.~66-71).
Here $\lambda _n $ runs over the linear envelope of the countable set
$\Lambda\subset {\bf R}$ over the field ${\bf Q}$.
The following properties of
Bochner-Fejer kernels are fulfilled:

1) $0 \le k_n^q  \le 1$;

2) there is only a finite number of nonzero coefficients $k_n^q $
for every fixed $q$;

3)  $k_n^q\to 1$ as $q \to \infty $ for every fixed $n$;

4) $K^{\left( q \right)} ( - t) = K^{\left( q \right)} (t)$;

5) $K^{\left( q \right)} (t) \ge 0$ for all $t \in {\bf R}$;

6) $M_t \left\{ {K^{(q)} (t)} \right\} = 1$.

Clearly,
$M_t \left\{ {e^{ - i\lambda t} e^{i\mu t} }
\right\} = 0$ for $\lambda  \ne \mu $.
Hence for any finite sum
$Q(x,y)$ of type (\ref{f1}) with all exponents
$\lambda _n \in \Lambda$,
we have
$$(Q * K^{\left( q \right)} )(x,y) = M_t \left\{ {Q(x + t,y)K^{\left(
q \right)} (t)} \right\}
= \sum\limits_{n = 1}^N {c_n (y)} k_n^q e^{i\lambda _n x}.$$
It follows from  condition 3) that this sum
converges to $Q(x,y)$ as $q \to \infty $
uniformly on the strip $\Pi _{\left[ {a,b} \right]} $.

Take a function $f \in W_{\left[ {a,b} \right]}^1 $.
Fix $\varepsilon  > 0$ and take any sum $Q(x,y)$ of type
 (\ref{f1}) such that
$d_{[a,b]}^{W^1 } (f,Q) < \varepsilon $.
Then choose  $q$ such that
$$\mathop {\sup }\limits_{(x,y) \in \Pi _{\left[ {a,b} \right]} } \left|
{Q(x,y) - (Q * K^{(q)} )(x,y)} \right| < \varepsilon. $$
Put $(f * K^{\left( q \right)} )(x,y) =
M_t \left\{ {f(x + t,y)K^{\left( q \right)} (t)}
\right\}$.
It follows from (\ref{f5}) that the mean value does not change
under shift on $x\in {\bf R}$. Since
$$M_t \left\{ {f(x + t,y)e^{ - i\lambda t} } \right\} = M_t \left\{ {f(x + t,y)e^{ - i\lambda (x
+ t)} } \right\}e^{i\lambda x},$$
we see that $\left\{ {f *
K^{(q)} } \right\}(x,y)$ is a finite sum
of type (\ref{f1}).
Next, we have
$$|(f * K^{(q)} )(x,y)| \le \mathop {\underline {\lim } }\limits_{T
\to \infty }
\frac{1}{{2T}}\int\limits_{ - T}^T {|f(x + t,y)|K^{(q)} (t)dt}$$
for all $x \in {\bf R},\;y \in K$.
By Fatou's lemma, we obtain
$$\frac{1}{{2X}}\int\limits_{ - X}^X {|(f * K^{(q)} )(x + \tau ,y)|dx
\le \mathop
{\underline {\lim } }\limits_{T \to \infty } } \frac{1}{{2T}}
\int\limits_{ - T}^T
{\frac{1}{{2X}}\int\limits_{ - X}^X {|f(x + t + \tau ,y)
|K^{(q)} (t)dt} } dx \le $$
$$\mathop {\overline {\lim } }\limits_{T \to \infty }
\frac{1}{{2T}}\int\limits_{ - T}^T
{K^{(q)} (t)} \left\{ {\mathop {\sup }\limits_{t \in {\bf R},
\;y \in \left[ {a,b} \right]}
\frac{1}{{2X}}\int\limits_{ - X}^X {|f(x + t,y)|dx} } \right\}dt$$
for every $\tau  \in {\bf R}$ and $X < \infty $.
The expression in the curly brackets does not exceed
$d_{[a,b]}^{W^1 } (f,0) + \varepsilon $ for large $X$.
Consequently, passing to the limit
as $X \to \infty $, and then
$\varepsilon  \to 0$, gives the inequality
$$d_{[a,b]}^{W^1 } (f * K^{(q)} ,0) \le d_{[a,b]}^{W^1 }
(f,0)\mathop {\lim }\limits_{T
\to \infty } \frac{1}{{2T}}\int\limits_{ - T}^T {K^{(q)} (t)dt}
= d_{[a,b]}^{W^1 } (f,0).$$
Replace $f$ by $f - Q$ to get
$$d_{[a,b]}^{W^1 } (f * K^{(q)} ,f) \le d_{[a,b]}^{W^1 }
((f - Q) * K^{(q)} ,0) +
d_{[a,b]}^{W^1 } (Q * K^{(q)} ,Q) + d_{[a,b]}^{W^1 }
(f - Q,0) \le 3\varepsilon,$$
thus the exponential sums $f * K^{(q)} $ approximate the function
$f$ with respect to the metric $d_{[a,b]}^{W^1 } $ as well.

{\bf The proof of Theorem 1.}

Let $f \in HW_{\left( {a,b} \right)}^1 $,
 $\left[
{\alpha ,\beta } \right] \subset (a,b)$.
According with Proposition 1 we have
 $\mathop {\sup }\limits_{z \in \Pi _{\left[ {\alpha
,\beta } \right]} } |f(z)| = C < \infty $,
therefore for any  $K^{(q)} $
$$\left| {\left( {f * K^{\left( q \right)} } \right)(z)} \right|
= \left| {M_t \left\{ {f(z +
t)K^{\left( q \right)} (t)} \right\}} \right| \le CM_t
\left\{ {K^{\left( q \right)} (t)} \right\} =
C.$$
Hence for any $p > 1$
$$\left( {d_{[\alpha ,\beta ]}^{W^p } (f * K^{(q)} ,f)} \right)^p
= \mathop {\overline {\lim
} }\limits_{T \to \infty } \mathop {\sup }\limits_{z \in
\Pi _{[\alpha ,\beta ]} }
\frac{1}{{2T}}\int\limits_{ - T}^T {|f(z + t) -
(f * K^{(q)} )(z + t)|^p dt} $$
$$ \le \mathop {\overline {\lim } }\limits_{T \to \infty }
 \mathop {\sup }\limits_{z \in \Pi
_{[\alpha ,\beta ]} } \frac{1}{{2T}}\int\limits_{ - T}^T {|f(z + t) - (f * K^{(q)} )(z + t)| \cdot
|f(z + t) - (f * K^{(q)} )(z + t)|^{p - 1} dt} $$
$$ \le \left( {2C} \right)^{p - 1} d_{[\alpha ,\beta ]}^{W^1 } (f * K^{(q)} ,f).$$
Thus the exponential sums $f * K^{(q)} $ approximate $f$ in the metric
$d_{[\alpha ,\beta ]}^{W^p } $, as claimed.

{\bf The proof of Theorem 2.}

Consider the function
 $f(z) = \sum\limits_{n \in I} {e^{ - 4(z - n)^2 } } $, where $I = \left\{
{n = 3^{l - 1} (3k + 1),\;k \in {\bf Z},\;l \in {\bf N}} \right\}$.
For any  $z = x + iy \in {\bf C}$ we have
$\left| {f(z)} \right| \le e^{4y^2 }
\sum\limits_{n \in I}^{} {e^{ - 4(x - n)^2 } }
 \le e^{4y^2 } \sum\limits_{n \in {\bf Z}}^{} {e^{ - 4(x
- n)^2 } } $.
The series
 $\sum\limits_{n =  - \infty }^\infty  {e^{ - 4(x - n)^2 } } $
converges for every $x \in {\bf R}$ and
is a periodic function  with the period 1, therefore $f(z)$ is bounded
in every strip $\left\{ {z = x + iy:x \in {\bf R},\;\left|
y \right| < H} \right\}$ and $f(z)$ is an entire function.
In particular, $f(z)$ is uniformly continuous on each strip.

Let us check that $f(z) \in HW_{( - \infty ,\infty )}^1 $.
Put
$$\varphi _l (z) = \sum\limits_{\scriptstyle n = 3^{l - 1}
(3k + 1), \atop   \scriptstyle k \in {\bf Z}}^{} {e^{ - 4(z - n)^2 } } $$.
The function $f_m (z) = \sum\limits_{1
\le l \le m}^{} {\varphi _l (z)} $
is the sum of the functions with periods $3^l ,\;l \le m$,
therefore $f_m(z)$ is
a periodic function with the real period $3^m $,
and  belongs to the space
 $HU_{\left( { - \infty ,\infty } \right)} $.
It is sufficient to prove that
\begin{equation}
\label{f6}
d_{\left[ { - H,H} \right]}^{W^p } \left( {f,f_m } \right) \to 0
\end{equation}
as $m \to \infty $ for each $H < \infty $.
We have
$$\left| {f(z) - f_m (z)} \right| = \sum\limits_{l = m + 1}^\infty
{\sum\limits_{n \in {\bf Z}}^{}
{\left| {e^{ - 4(z - 3^l n - 3^{l - 1} )^2 } } \right|} }  \le \sum\limits_{n =  - \infty }^\infty  {e^{
- 4(x - 3^m n)^2 } e^{4H^2 } } $$
and
$d_{\left[ { - H,H} \right]}^{W^1 } \left( {f,f_m } \right)
 \le e^{4H^2 } \overline {\mathop
{\lim }\limits_{T \to \infty } } \mathop
{\sup }\limits_{x \in {\bf R}} \frac{1}{{2T}}\sum\limits_{n =
- \infty }^\infty  {\int\limits_{ - T}^T
{e^{ - 4(x + t - 3^m n)^2 } dt} } $.

Put $$
E_1 = \left\{ {n \in {\bf Z}:
n \le 3^{ - m} x - 3^{ - m} T - \frac{1}{2}} \right\}, \
n_1  = \mathop {\sup }\limits_{} E_1, $$
$$E_2  = \left\{ {n \in {\bf Z} : 3^{ - m} x - 3^{ - m} T - \frac{1}{2} < n < 3^{ - m} x + 3^{ - m} T
+ \frac{1}{2}} \right\},$$
$$E_3  = \left\{ {n \in {\bf Z} :n \ge 3^{ - m} x + 3^{ - m} T + \frac{1}{2}} \right\},
\ n_2  = \mathop {\inf }\limits_{} E_3$$
for any fixed $x \in {\bf R},\;\;T > 0.$
Denote by ${\rm card}\,E$ the number  of elements
of the set $E$.
Note that ${\rm card}\,E_2 \le 2\cdot 3^{ - m} T + 2$.

For  $n \in E_1 $ and $t \in \left[ { - T,T} \right]$
$$\left( {x + t - 3^m n} \right) \ge 3^m \left( {n_1  + \frac{1}{2} - n} \right)$$
and
$$e^{ - 4(x + t - 3^m n)^2 }  \le e^{ - 3^{2m} (2(n_1  - n) + 1)^2 }. $$
Similarly,
$$\left( {3^m n - x - t} \right) \ge 3^m \left( {n - n_2  + \frac{1}{2}} \right)$$
and
$$e^{ - 4(x + t - 3^m n)^2 }  \le e^{ - 3^{2m} (2(n - n_2 ) + 1)^2 }$$
for $n \in E_3 $, $t \in \left[ { - T,T} \right]$.
Consequently,
$$\sum\limits_{n \in E_1 }^{} {\frac{1}{{2T}}\int\limits_{ - T}^T {e^{ - 4(x + t - 3^m
n)^2 } dt}  \le \sum\limits_{n = 1}^\infty  {e^{ - (3^m n)^2 } }  \le \int\limits_0^\infty  {e^{ -
(3^m u)^2 } du}  = \frac{{\sqrt \pi  }}{{2 \cdot 3^m }}},$$
and the same estimate holds for the sum
$$\sum\limits_{n \in E_3 }^{} {\frac{1}{{2T}}\int\limits_{ - T}^T {e^{ - 4(x + t - 3^m
n)^2 } dt} }. $$
Next, we have
$$\int\limits_{ - T}^T {e^{ - 4(x + t - 3^m n)^2 } dt}  \le \int\limits_{ - \infty }^\infty  {e^{
- 4u^2 } du}  = \frac{{\sqrt \pi  }}{2}$$
and
$$\sum\limits_{n \in E_2 }^{} {\frac{1}{{2T}}\int\limits_{ - T}^T
{e^{ - 4(x + t - 3^m
n)^2 } dt}  \le \frac{1}{{2T}} \cdot \frac{{\sqrt \pi
 }}{2}{\rm card}\,E_2  \le \frac{1}{{2T}}
\cdot \frac{{\sqrt \pi  }}{2}\left( {3^{ - m} 2T + 2} \right)} $$
for $n \in E_2 $.
Thus $d_{[ - H,H]}^{W^1 } \left( {f,f_m } \right) \le
\frac{{3\sqrt \pi  }}{2}
\cdot 3^{ - m} e^{4H^2 } $, and we obtain (\ref{f6}).

Choose $H < \infty $.
Let us check that $d_{[ - H,H]}^{W^1 } (f,g) = 0$
for no function
$g \in HU_{\left( { - \infty ,\infty }
\right)} $.
We will prove the stronger result:
there are no almost periodic functions $g(x)$ in the sense of Stepanov
 on ${\bf R}$ with the property
 $d_{\left\{ 0 \right\}}^{W^1 } (f,g) = 0$.

We need some auxiliary lemmas:

{\bf Lemma 1.}

\begin{equation}
\label{f7}
\mathop {\sup }\limits_{x \in {\bf Z}\backslash I} f(x)
\le \frac{{\sqrt \pi  }}{2} < 1 \le \mathop
{\inf }\limits_{x \in I} f(x).
\end{equation}

{\bf The proof of Lemma 1.}

We have
$$f(x) = \sum\limits_{l \in {\bf N}} {\varphi _l } (z) \le
\sum\limits_{n \in {\bf Z}\backslash \{ 0\} }^{} {e^{ - 4n^2 } }
\le \int\limits_{ - \infty }^\infty  {e^{
- 4t^2 } dt}  = \frac{{\sqrt \pi  }}{2}$$
for any $x \in {\bf Z}\backslash I$,
and $
f(x)\ge e^{-4(x-n_0)^2}=1$  for $x = n_0\in I$.

{\bf Lemma 2.}

{\it For all $q \in {\bf Z} \backslash \left\{ 0 \right\}$
there exists a two-sided arithmetical progression
 $I\left( q \right) \subset I$

such that $\left( {I(q) + q} \right)
\cap I = \emptyset $.}

{\bf The proof of Lemma 2.}

Every positive integer $q$ has the form $q = 3^{r - 1} (3m + 1)$ or
$q = 3^{r - 1} (3m - 1)$, $r \in {\bf N},\;m \in {\bf Z}$.
In the first case take
 $n_j  = 3^{r - 1} (3j + 1)$, $j \in {\bf Z}$,
because
$$3^{r - 1} (3m + 1) + 3^{r - 1} (3j + 1) =
3^{r - 1} \left( {3(m + j + 1) - 1} \right) \notin I.$$
In the second case take
$n_j  = 3^{r - 1} (3(3j - m - 1) + 1)$, $j \in {\bf Z}$, because
$n_j  + q = 3^r (3j - 1) \notin I$ for any $j \in {\bf Z}$.

{\bf Lemma 3.}

{\it There exists $\gamma  > 0$ such that
for each $\tau  \in {\bf R}$, $\left| \tau
\right| \ge 1$ there is a
relatively dense set $I(\tau ) \subset {\bf R}$,
with the property
\begin{equation}
\label{f8}
\inf_{x\in I(\tau)}\left| {f(x + \tau ) - f(x)} \right| > \gamma.
\end{equation}}
{\bf The proof of Lemma 3.}

Since $f(x)$ is uniformly continuous
on ${\bf R}$, there exists $N < \infty $ such that
\begin{equation}
\label{f9}
|f(x) - f(t)| < \frac{1}{2}\left( {1 - \frac{{\sqrt \pi  }}{2}} \right)
\end{equation}
whenever $\left| {x - t} \right| \le \frac{1}{N}$.

Let $\tau $ be an arbitrary real number,
$\left| \tau  \right| \ge 1$.
We show that inequality (\ref{f8}) takes place
with
$\gamma  = \frac{1}{{2N}}\left( {1 - \frac{{\sqrt \pi  }}{2}} \right)$
for all points from some relatively dense set in ${\bf R}$.

Since the fractional parts of numbers
$0,\tau ,\;2\tau ,\;...,\;N\tau $
belong to the half-open interval
 $[0,1)$, there are two numbers
$k\tau $ and $k'\tau $, $0 \le k < k' \le N$,
such that the distance between their  fractional parts
is at most $\frac{1}{N}$, i.e.
for some
$q \in {\bf Z}\backslash \left\{ 0 \right\}$
the inequality
$$\left| {\,k\tau  - k'\tau  - q} \right| \le \frac{1}{N}$$
holds.
For $M = k' - k$ we obtain
\begin{equation}
\label{f10}
\left| {\,M\tau  - q} \right| \le \frac{1}{N}.
\end{equation}

Let $L$ be the difference of
the arithmetic progression $I(q)$ from Lemma 2. Fix a real number
 $a \in {\bf R}$, and $n \in I(q) \cap [a,a + L)$.
Taking into account Lemma 1 and 2, we see that
\begin{equation}
\label{f11}
\left| {f(n) - f(n + q)} \right| \ge 1 - \frac{{\sqrt \pi  }}{2}.
\end{equation}

On the other hand,
\begin{equation}
\label{f12}
\left| {f(n + q) - f(n)} \right| \le |f(n + q) - f(n + M\tau )| + \sum\limits_{k = 0}^{M - 1}
{\left| {f(n + k\tau ) - f(n + (k + 1)\tau )} \right|}.
\end{equation}

By (\ref{f9}) and (\ref{f10}), we have
$\left| {f(n + q) - f(n + M\tau )} \right| \le
\frac{1}{2}\left( {1 - \frac{{\sqrt \pi  }}{2}} \right)$.
Hence inequalities (\ref{f11}) and (\ref{f12})
imply that
$$|f(n + k''\tau ) - f(n + (k'' + 1)\tau )|
\ge \frac{1}{{2N}}\left( {1 - \frac{{\sqrt \pi  }}{2}} \right)$$
 for some $k''$, $0 \le k'' \le M \le N$.

Thus there exists the point $x = n + k''\tau $,  $x\in [a,a + L + N\tau ]$
 such that
$\left| {f(x + \tau ) - f(x)} \right| > \gamma$. The lemma is proved.

We continue the proof of Theorem 3. Let $\gamma$ be the constant from Lemma 3.
Take $\delta  > 0$ such that the inequality
\begin{equation}
\label{f13}
\left| {f(x) - f(t)} \right| \le \frac{\gamma }{5},
\end{equation}
holds whenever $x,\;t \in {\bf R}$, $\left| {x - t} \right|<\delta$.
Let us check that an arbitrary function $g(x)$ from the equivalent class
of the function $f(x)$ in the space $W_{\{ 0\} }^1 $
satisfies the inequality
\begin{equation}
\label{f14}
d_{\{ 0\} }^{S^1 } \left( {g(t + \tau ),g(t)} \right)
= \mathop {\sup }\limits_{x \in {\bf R}}
\int\limits_x^{x + 1} {|g(t + \tau ) - g(t)|dt}
\ge \frac{{\gamma \delta }}{{10}}
\end{equation}
for arbitrary $\tau  \in {\bf R}$, $\left| \tau  \right| \ge 1$.
Then the set of $\varepsilon$-almost periods for the function $g$
in the Stepanov metric for
$\varepsilon <\frac{{\gamma \delta }}{{10}}$ is contained in the segment
$[-1,1]$, and so $g(x)\notin S^1_{\{0\}}$.

In order to prove (\ref{f14}) for fixed $\tau  \in {\bf R}$
put
$$F_1  = \left\{ {x \in {\bf R}:\left| {g(x) - f(x)} \right| \ge \frac{\gamma }{5}} \right\},$$
$$F_2  = \left\{ {x \in {\bf R}:\left| {g(x + \tau ) - f(x + \tau )} \right| \ge \frac{\gamma }{5}}
\right\}.$$
Take $L < \infty $ such that the set $I(\tau)$ from Lemma 3 has
nonempty intersection
with every interval of the length $L$.
Since
$$\frac{1}{{2nL}}mes\left( {F_1  \cap \left[ { - nL,\;nL} \right]} \right) \le
\frac{5}{{2nL\gamma }}\int\limits_{F_1  \cap \left[ { - nL,\;nL} \right]}^{} {\left| {f(x) -
g(x)} \right|dx \le } \frac{5}{{2nL\gamma }}\int\limits_{ - nL}^{nL} {\left| {f(x) - g(x)}
\right|dx}, $$
and $f$, $g$ belong to the same equivalent class in the
space $W^1_{\{0\}}$, we get
$$\mathop {\overline {\lim } }\limits_{n \to \infty }
\frac{1}{{2nL}}mes\left( {F_1  \cap
\left[ { - nL,\;nL} \right]} \right) \le d_{\{ 0\} }^{W^1 } \left\{ {f,g} \right\} = 0.$$
The same equality holds for the set
$F_2 $. Hence
$$\mathop {\lim }\limits_{n \to \infty } \frac{1}{{2n}}\sum\limits_{k =  - n}^{n - 1}
{\frac{{mes\left\{ {\left( {F_1  \cup F_2 } \right) \cap [kL,(k + 1)L]} \right\}}}{L}}  =
\mathop {\lim }\limits_{n \to \infty } \frac{1}{{2nL}}mes\left\{ {\left( {F_1  \cup F_2 } \right)
\cap [ - nL,nL]} \right\} = 0.$$
Therefore for  $n$ sufficiently large
there exists an interval $[k_0 L,(k_0  + 1)L]$ such that
\begin{equation}
\label{f15}
mes\left\{ {\left( {F_1  \cup F_2 } \right) \cap
[k_0 L,(k_0  + 1)L]} \right\} < \frac{\delta
}{2}.
\end{equation}

Take a real number  $x \in [k_0 L,(k_0  + 1)L]\cap I(\tau)$,
where $I(\tau)$ is defined in Lemma 3.
 By (\ref{f13}) we get the inequality
$$\left| {f(t + \tau ) - f(t)} \right| \ge |f(x + \tau ) - f(x)|
- |f(x + \tau ) - f(t + \tau )| - |f(x) -
f(t)| \ge \frac{{3\gamma }}{5}$$
 for each point $t \in (x - \delta ,x + \delta)$.

Note that the length of the interval
 $(x - \delta ,x + \delta ) \cap (k_0 L,(k_0  + 1)L)$ is at least
$\delta $.
It follows from (\ref{f15}) that the measure of the set
 $(x - \delta ,x + \delta )\backslash \left[ {F_1
\cup F_2 } \right]$ is at least $\frac{\delta }{2}$.
Next, for $t\in (x - \delta ,x + \delta )\backslash \left[ {F_1
\cup F_2 } \right]$  we have
$$\left| {g(t + \tau ) - g(t)} \right| \ge \left| {f(t + \tau ) - f(t)} \right| - \left| {f(t + \tau ) - g(t
+ \tau )} \right| - \left| {f(t) - g(t)} \right| \ge \frac{\gamma }{5}.$$
Thus,
$$\int\limits_{x - \delta }^{x + \delta } {\left| {g(t + \tau ) - g(t)} \right|dt \ge }
\int\limits_{[x - \delta ,x + \delta ]\backslash (F_1  \cup F_2 )}^{} {\left| {g(t + \tau ) - g(t)}
\right|dt \ge } \frac{{\gamma \delta }}{{10}}.$$
The last inequality implies (\ref{f14}). The theorem is proved.

For the proof of other theorems we need following lemmas:

{\bf Lemma 4.}

{\it Any collection of functions
 $\ae _n (x) = e^{ - 4(x - 3n)^2 } $ satisfy the inequality
$$ \left({\sum\limits_{k = 1}^\infty  {\ae _k (x)} }\right)^p \le 2^{p-1}
\left( {\sum\limits_{k = 1}^\infty  {{\ae _k (x)} } } \right).$$}

{\bf The proof of Lemma 4.}

Fix $x \in {\bf R}$ and put
$n_0  = \left[ {\frac{x}{3} + \frac{1}{2}}
\right]$. Since $3n_0  -
\frac{3}{2} \le x \le 3n_0  - \frac{1}{2}$, we obtain
$$\sum\limits_{n =  - \infty }^{n_0 } {e^{ - 4(x - 3n)^2 } }
\le \sum\limits_{n =  - \infty
}^{n_0 } {e^{ - (6(n_0 - 1  - n) + 3)^2 } }  \le
\sum\limits_{n = 1}^\infty  {e^{ - (3n)^2 } }  \le
\frac{{\sqrt \pi  }}{6}$$
and
$$\sum\limits_{n = n_0+1 }^\infty  {e^{ - 4(x - 3n)^2 } }
\le \sum\limits_{n = n_0+1 }^\infty
{e^{ - (6(n - n_0 -1 ) + 3)^2 } }  \le \sum\limits_{n = 1}^\infty  {e^{ - (3n)^2 } }  \le \frac{{\sqrt
\pi  }}{6},$$
hence
 $\sum\limits_{n \in {\bf Z}\backslash \left\{ {n_0}
\right\}}^\infty  {e^{ - 4(x - 3n)^2 }
}  \le \frac{{\sqrt \pi  }}{3} < 1$.
The statement of the lemma follows from the inequality
$$
(a+b)^p\le 2^{p-1}(a+b)
$$
with $a=\ae _{n_0} (x)$ and
$b=\sum\limits_{n \neq n_0} \ae _n (x).
$

{\bf The proof of Theorem 3.}

For $z = x + iy \in {\bf C},\;l \in {\bf N}$ put
$$\varphi _l (z) = \sum\limits_{\scriptstyle n = 3^{l - 1} (3k + 1), \atop
  \scriptstyle k \in {\bf Z}} {e^{ - 4(z - n)^2 } },$$
\begin{equation}
\label{f16}
f(z) = \sum\limits_{l = 1}^\infty  {l\varphi _l (z)}.
\end{equation}

First of all, for $\left| y \right| \le H$ and any $x \in {\bf R}$
we have
\begin{equation}
\label{f17}
\left| {\varphi _l (z)} \right| \le e^{4H^2 }
\sum\limits_{k \in {\bf Z}} {e^{ - 4(x - 3^{l - 1} (3k +
1))^2 } }  \le e^{4H^2 }
\sum\limits_{n \in {\bf Z}}^{} {e^{ - 4(x - n)^2 } }.
\end{equation}
Hence each term of sum (\ref{f16}) is uniformly bounded
 on every horizontal strip.
If $|x| \le 3^{l_0  - 2} $ for some $l_0  \in {\bf N}$,
 then for all $k \in {\bf Z},\;l \ge l_0 $ we have
$$\left| {x - 3^{l - 1} (3k + 1)} \right| \ge 3^{l - 1} |3k + 1| - 3^{l - 2}  \ge 3^{l - 2} (4|k| +
1)$$
and for $z = x + iy$, $|x| \le 3^{l_0  - 2} $, $|y| \le H$
$$\left| {\varphi _l (z)} \right| \le e^{4H^2 } \sum\limits_{k \in {\bf Z}}^{} {e^{ - 4(3^{l - 2} (4|k|
+ 1))^2 } }  \le e^{4H^2 } 2\sum\limits_{n = 1}^\infty  {e^{ - 4(3^{l - 2} n)^2 } }  \le e^{4H^2
} 2\int\limits_0^\infty  {e^{ - 4(3^{l - 2} u)^2 } du}  = \frac{{e^{4H^2 } \sqrt \pi  }}{{2 \cdot
3^{l - 2} }},$$
so all terms of series (\ref{f16})
with indices $l \ge l_0 $ are majorized by the terms of  the
convergent series
$$\sum\limits_{l = 1}^\infty  {\frac{{9\sqrt \pi
e^{4H^2 } }}{2} \cdot \frac{l}{{3^l }}}.
$$
Thus series (\ref{f16}) uniformly converge
on every compact set in ${\bf C}$, and
$f(z)$ is an entire function.

Next,  $\varphi _l (z)$ is an entire function with the period
$3^l $, and the sum
$\sum\limits_{l = 1}^m l \varphi _l (z)$
is a periodic function with the period $3^m $,
therefore this sum belongs to the space
$HU_{\left( { - \infty,\infty } \right)} $.
Hence  if
\begin{equation}
\label{f18}
\overline {\mathop {\lim }\limits_{T \to \infty } } \left( {\frac{1}{{2T}}\mathop {\sup
}\limits_{|y| \le H} \int\limits_{ - T}^T {\left| {\sum\limits_{l = m + 1}^\infty  {l\varphi _l (z)}
dx} \right|^p } } \right)^{\frac{1}{p}}
= d_{[ - H,H]}^{B^p } \left( {\sum\limits_{l = 1}^m {l\varphi _l (z),f(z)} }
 \right) \to 0
\end{equation}
as $m \to \infty $ for all $H > 0$, then
$f(z)\in HB^p_{(-\infty, \infty)}$, $p \ge 1$.

Fix $T > 0$ and consider the integral
$$\int\limits_{ - T}^T {\varphi _l (x)dx = \sum\limits_{k \in {\bf Z}}^{} {\int\limits_{ - T}^T
{e^{ - 4(x - 3^{l - 1} (3k + 1))^2 } dx} } }.$$
 Put
$$E_1  = \left\{ {n \in {\bf Z}:n \le  - 3^{1 - l} T - \frac{1}{2}} \right\},
\ n_1  = \mathop {\sup }\limits_{} E_1,$$
$$E_2  = \left\{ {n \in {\bf Z}: - 3^{1 - l} T - \frac{1}{2} < n < 3^{1 - l} T + \frac{1}{2}}
\right\}\backslash \left\{ 0 \right\},$$
$$E_3  = \left\{ {n \in {\bf Z}:n \ge 3^{1 - l} T + \frac{1}{2}} \right\},
\ n_2  = \mathop {\inf }\limits_{} E_3.$$

Note that
$$\sum\limits_{k \in {\bf Z}}^{} {\int\limits_{ - T}^T {e^{ - 4(x - 3^{l - 1}
 (3k + 1))^2 } dx} }
\le \sum\limits_{n \in E_1 }^{} {\int\limits_{ - T}^T {e^{ - 4(x - 3^{l - 1}
n)^2 } dx} }  + $$
$$
\sum\limits_{n \in E_2 }^{} {\int\limits_{ - T}^T {e^{ - 4(x - 3^{l - 1} n)^2
 } dx} }  +
\sum\limits_{n \in E_3 }^{} {\int\limits_{ - T}^T {e^{ - 4(x - 3^{l - 1} n)^2
 } dx} }.$$
Note also that the number $card\,E_2 $ is at most  $2T \cdot 3^{1 - l} $.

For $n \in E_1 $ and $x \in \left[ { - T,T} \right]$
we have
$$\left( {x - 3^{l - 1} n} \right) \ge 3^{l - 1} \left( {n_1  + \frac{1}{2} - n} \right)$$
and
$$e^{ - 4(x - 3^{l - 1} n)^2 }  \le e^{ - 3^{2(l - 1)}
(2(n_1  - n) + 1)^2 }, $$
by the same way, for $n \in E_3 $, $x \in \left[ { - T,T} \right]$,
we have
$$\left( {3^{l - 1} n - x} \right) \ge 3^{l - 1} \left( {n - n_2  + \frac{1}{2}} \right)$$
and
$$e^{ - 4(x - 3^{l - 1} n)^2 }  \le e^{ - 3^{2(l - 1)} (2(n - n_2 ) + 1)^2 }. $$
Therefore
$$\sum\limits_{n \in E_1 }^{} {\frac{1}{{2T}}
\int\limits_{ - T}^T {e^{ - 4(x - 3^{l - 1}
n)^2 } dx}  \le \sum\limits_{n = 1}^\infty
{e^{ - 3^{2(l - 1)} n^2 } }  \le \int\limits_0^\infty
{e^{ - (3^{l - 1} u)^2 } du}  =
\frac{{\sqrt \pi  }}{{2 \cdot 3^{l - 1} }}}.$$
The same estimate is true for the sum over $n \in E_3 $.
Next, for $n \in E_2 $ we get
$$\int\limits_{ - T}^T {e^{ - 4(x - 3^{l - 1} n)^2 } dx}  \le \int\limits_{ - \infty }^\infty
{e^{ - 4u^2 } du}  = \frac{{\sqrt \pi  }}{2}$$
and
$$\sum\limits_{n \in E_2 }^{}
{\frac{1}{{2T}}\int\limits_{ - T}^T {e^{ - 4(x - 3^{l - 1}
n)^2 } dx}
\le \frac{{\sqrt \pi  }}{{4T}}{\rm card}\,E_2
\le \frac{{\sqrt \pi  }}{{4T}}3^{1 - l} 2T.} $$
Thus $\frac{1}{{2T}}\int\limits_{ - T}^T {\varphi _l (x)dx =
\sum\limits_{k
\in {\bf Z}}^{} {\frac{1}{{2T}}\int\limits_{ - T}^T {e^{ - 4(x - 3^{l - 1} (3k + 1))^2 } dx}  \le
\frac{{3\sqrt \pi   \cdot 3^{1 - l} }}{2}} } $.

Applying Lemma 4 to the functions
 $e^{ - 4(x - 3^{l - 1} (3k + 1))^2 } $, $k \in
{\bf Z}$,
we obtain
\begin{equation}
\label{f19}
\frac{1}{{2T}}\int\limits_{ - T}^T {\varphi _{_l }^p (x)dx \le 2^p
\frac{1}{{2T}}\int\limits_{ - T}^T {\varphi _{_l }^{} (x)dx}  \le 2^{p - 1} } 9\sqrt \pi  3^{ -
l}.
\end{equation}
The H\"older inequality implies
$$\left| {\sum\limits_{l = m + 1}^\infty  {l\varphi _l (x)} } \right| \le \left( {\sum\limits_{l
= m + 1}^\infty  {l^{2p} \varphi _l^p (x)} } \right)^{\frac{1}{p}} \left( {\sum\limits_{l = m +
1}^\infty  {\frac{1}{{l^q }}} } \right)^{\frac{1}{q}},$$
where $\frac{1}{p} + \frac{1}{q} = 1.$
Therefore for sufficiently large $m$ we have
$$\frac{1}{{2T}}\int\limits_{ - T}^T {\left( {\sum\limits_{l = m + 1}^\infty  {l\varphi _l
(x)} } \right)} ^p dx \le \sum\limits_{l = m + 1}^\infty  {l^{2p} \frac{1}{{2T}}\int\limits_{ -
T}^T {\varphi _{_l }^p (x)} } dx \le 2^{p - 1} 9\sqrt \pi  \sum\limits_{l = m + 1}^\infty
{\frac{{l^{2p} }}{{3^l }}}.$$
If we combine the letter with the inequality
$$\mathop {\sup }\limits_{|y| \le H} \left|
{\sum\limits_{l = m + 1}^\infty  {l\varphi _l (z)}
} \right| \le e^{4H^2 } \sum\limits_{l = m + 1}^\infty
{l\varphi _l (x)},$$
we get (\ref{f18}).

Let us check that an arbitrary function $g(x)$ from the equivalent class
of $f(x)$ in the space $B_{\{ 0\} }^p $
does not belong to the space $W_{\{ 0\} }^1 $.
Then it follows that $g$ does not belong to the space $HW^1_{[-H,H]}$
for all $H>0$.

Assume the contrary. Then
$$\mathop {\overline {\lim } }\limits_{x \to \infty } \frac{1}{{2T}}\mathop {\sup
}\limits_{x \in {\bf R}} \int\limits_{ - T}^T {|g(x + t)|^p dt = \left( {d_{\left\{ 0 \right\}}^{W^p }
(0,g)} \right)^p  < \infty }, $$
and for some
 $T_0  < \infty $, $c = c(T_0 )$ and all $x \in {\bf R}$
\begin{equation}
\label{f20}
\int\limits_{ - T_0 }^{T_0 } {|g(x + t)|^p dt \le c}.
\end{equation}

Take an integer $l$ such that
\begin{equation}
\label{l}
l > \mathop {\max }\limits_{} \left\{ {\left(
{\frac{{2c}}{{\int\limits_{ - T_0 }^{T_0 }
{e^{ - 4pt^2 } dt} }}} \right)^{\frac{1}{p}}
,\frac{{\log 2T_0 }}{{\log 3}}} \right\}.
\end{equation}

For each fixed
 $x_n  = 3^l n + 3^{l - 1} ,\;n \in {\bf Z}$,
we have
\begin{equation}
\label{f21}
\int\limits_{ - T_0 }^{T_0 } {|f(x_n  + t)|^p dt}
\ge l^p \int\limits_{ - T_0 }^{T_0 } {\left(
{{e^{ - 4(x_n  + t - 3^{l - 1} (3n + 1))^2 } } } \right)^p dt \ge l^p
\int\limits_{ - T_0 }^{T_0 } {e^{ - 4pt^2 } dt} }  \ge 2c.
\end{equation}

Take $T_n  = x_n  + T_0 $, $n \in {\bf N}$.
It follows from ($\ref{l}$) that $2T_0  < 3^l$,
hence the intervals $[x_k  - T_0 ,\;x_k  + T_0 ]$
are mutually disjoint.
By the Minkowski inequality for
 $E = \bigcup\limits_{k =  - n}^n {[x_k  - T_0 ,x_k
+ T_0 ]} $,
 (\ref{f20}) and (\ref{f21}), we get
$$\left( {\int\limits_{ - T_n }^{T_n } {\left| {f(t) - g(t)} \right|^p dt} }
\right)^{\frac{1}{p}}  \ge \left( {\int\limits_E^{} {\left| {f(t) - g(t)} \right|^p dt} }
\right)^{\frac{1}{p}}  \ge $$
$$\left( {\int\limits_E^{} {\left| {f(t)} \right|^p dt} } \right)^{\frac{1}{p}}  - \left(
{\int\limits_E^{} {\left| {g(t)} \right|^p dt} } \right)^{\frac{1}{p}}  \ge $$
$$\left( {\sum\limits_{k =  - n}^n {\int\limits_{ - T_0 }^{T_0 } {\left| {f(x_k  + t)}
\right|^p dt} } } \right)^{\frac{1}{p}}  - \left( {\sum\limits_{k =  - n}^n {\int\limits_{ - T_0
}^{T_0 } {\left| {g(x_k  + t)} \right|^p dt} } } \right)^{\frac{1}{p}}  \ge $$
$$\left( {(2n + 1)\,2c} \right)^{\frac{1}{p}}  - \left( {(2n + 1)\,c} \right)^{\frac{1}{p}}  =
(2n + 1)^{\frac{1}{p}} c^{\frac{1}{p}} (2^{\frac{1}{p}}  - 1),$$
therefore
\begin{equation}
\label{f22}
\int\limits_{ - T_n }^{T_n } {\left| {f(t) - g(t)} \right|^p dt}  \ge (2n + 1)c(2^{\frac{1}{p}}  -
1)^p.
\end{equation}

Since
 $T_n  = T_0  + 3^{l - 1}  + n3^{l - 1} $,
we see that $\frac{{2n + 1}}{{2T_n }} \to 3^{1 -
l} $ as $n \to \infty $.
Hence inequality (\ref{f22}) contradicts to the equality
$$\mathop {\overline {\lim } }\limits_{T \to \infty }
\frac{1}{{2T}}\int\limits_{ - T}^T
{\left| {f(t) - g(t)} \right|^p dt}  =
\left( {d_{\left\{ 0 \right\}}^{B^p } (f,g)} \right)^p  = 0,$$
which is true for each function $g$ from the equivalent class of $f$.
The theorem is proved.

{\bf The proof of Theorem 4.}

Fix $p_0  \in (p,p')$ and take
\begin{equation}
\label{f23}
f(z) = \sum\limits_{l = 1}^\infty  {3^{\frac{l}{p_{0}}} \varphi_l (z)},
\end{equation}
here the functions $\varphi _l (z)$ are the same as in the proof
of Theorem 3. As in that proof we see that the terms of
series (\ref{f23}) are majorized by the terms of series
$$\sum\limits_{l = 1}^\infty  {\frac{{9\sqrt \pi  e^{4H^2 } 3^{l\left( {\frac{1}{{p_0 }} -
1} \right)} }}{2}} $$
on compact sets
$\left\{ {\left| x \right| \le 3^{l_0  - 2} ,
\;\left| y \right| \le H} \right\}$.
Therefore $f(z)$ is an entire function on  ${\bf C}$.

As in the proof of Theorem 3, we have for any  $H < \infty $
\begin{equation}
\label{f24}
d_{[ - H,H]}^{B^p } \left( {\sum\limits_{l = 1}^m {3^{\frac{l}{{p_0 }}} \varphi _l
(z),f(z)} } \right) \to 0\
 as \ m \to \infty .
\end{equation}

Indeed, by the H\"older equality
$$\left| {\sum\limits_{l = m + 1}^\infty  {3^{\frac{l}{{p_0 }}} \varphi _l (x)} } \right| \le
\left( {\sum\limits_{l = m + 1}^\infty  {3^{\frac{l}{{p_0 }}} l^p \varphi _l^p (x)} }
\right)^{\frac{1}{p}} \left( {\sum\limits_{l = m + 1}^\infty  {\frac{1}{{l^q }}} }
\right)^{\frac{1}{q}}, \ \frac{1}{p} + \frac{1}{q} = 1.$$

Hence for $m$ sufficiently large, for all $T > 0$
inequality (\ref{f19}) implies
$$\frac{1}{{2T}}\int\limits_{ - T}^T {\left( {\sum\limits_{l = m + 1}^\infty
{3^{\frac{l}{{p_0 }}} \varphi _l (x)} } \right)} ^p dx \le \sum\limits_{l = m + 1}^\infty
{3^{\frac{l}{{p_0 }}} l^p \frac{1}{{2T}}\int\limits_{ - T}^T {\varphi _{_l }^p (x)} } dx \le
2^{p - 1} 9\sqrt \pi  \sum\limits_{l = m + 1}^\infty  {\frac{{l^p }}{{3^{l\left( {1 -
\frac{1}{{p_0 }}} \right)} }}}.$$

The convergence of the last series yields
(\ref{f24}).

Let us show that an arbitrary function
 $g$ from the equivalent class
of  $f(x)$ in the space $B_{\{ 0\} }^p $
does not belong to $ B_{\left\{ 0 \right\}}^{^{p'} }$ .

Assume the contrary. Then for sufficiently large $T$ we have
$$\frac{1}{{2T}}\int\limits_{ - T}^T {|g(t)|^{p'} dt}  \le C_1  < \infty.$$
Let $x_n $ be the same as in the proof of Theorem 3, and
$T_n  = x_n  + \frac{1}{2}$.

Applying the H\"older inequality for the set
$${ E = }\bigcup\limits_{k =  - n}^n {\left[ {x_k  - \frac{1}{2},\;x_k  + \frac{1}{2}}
\right]}, $$
gives for sufficiently large  $n$
$$\int\limits_E^{} {|g(t)|^p dt \le \left( {\int\limits_E^{} {|g(t)|^{p'} dt} }
\right)^{\frac{p}{{p'}}} \left( {\int\limits_E^{} {dt} } \right)^{1 - \frac{p}{{p'}}}  \le } $$
$$(2T_n )^{\frac{p}{{p'}}} \left( {\frac{1}{{2T_n }}\int\limits_{ - T_n }^{T_n }
{|g(t)|^{p'} dt} } \right)^{\frac{p}{{p'}}} \left( {\int\limits_E^{} {dt} } \right)^{1 -
\frac{p}{{p'}}}  \le $$
$$C_1^{^{\frac{p}{{p'}}} } \left( {2T_n } \right)^{\frac{p}{{p'}}} (2n + 1)^{1 -
\frac{p}{{p'}}}. $$

As in the previous theorem, we have for all $l \in {\bf N}$
$$\int\limits_{ - \frac{1}{2}}^{\frac{1}{2}} {|f(x_n  + t)|^p dt}
\ge 3^{\frac{{pl}}{{p_0
}}} \int\limits_{ - \frac{1}{2}}^{\frac{1}{2}} {\left(
{\sum\limits_{k \in {\bf Z}}^{} {e^{ - 4(x_n  +
t - 3^{l - 1} (3k + 1))^2 } } } \right)^p dt
\ge 3^{\frac{{lp}}{{p_0 }}} \int\limits_{ -
\frac{1}{2}}^{\frac{1}{2}} {e^{ - 4pt^2 } dt} }  =
 3^{\frac{{lp}}{{p_0 }}} C_2.$$
Now the Minkowski inequality yields
$$\left( {\int\limits_{ - T_n }^{T_n } {\left| {f(t) - g(t)} \right|^p dt} }
\right)^{\frac{1}{p}}  \ge \left( {\int\limits_E^{} {\left| {f(t) - g(t)} \right|^p dt} }
\right)^{\frac{1}{p}}  \ge $$
$$\left( {\int\limits_E^{} {\left| {f(t)} \right|^p dt} } \right)^{\frac{1}{p}}  - \left(
{\int\limits_E^{} {\left| {g(t)} \right|^p dt} } \right)^{\frac{1}{p}}  \ge $$
$$\left( {\sum\limits_{k =  - n}^n {\int\limits_{ - \frac{1}{2}}^{\frac{1}{2}} {\left| {f(x_k
+ t)} \right|^p dt} } } \right)^{\frac{1}{p}}  - \left( {\int\limits_E^{} {\left| {g(t)} \right|^p dt}
} \right)^{\frac{1}{p}}  \ge $$
$$\left( {3^{\frac{{lp}}{{p_0 }}} C_2 (2n + 1)} \right)^{\frac{1}{p}}  - \left(
{C_1^{^{\frac{p}{{p'}}} } \left( {\frac{{2T_n }}{{2n + 1}}} \right)^{\frac{p}{{p'}}} (2n +
1)} \right)^{\frac{1}{p}}, $$
therefore
$$\left( {\frac{1}{{2T_n }}\int\limits_{ - T_n }^{T_n } {\left| {f(t) - g(t)} \right|^p dt} }
\right)^{\frac{1}{p}}  \ge \left( {\frac{{3^{\frac{{lp}}{{p_0 }}} C_2 (2n + 1)}}{{2T_n }}}
\right)^{\frac{1}{p}}  - \left( {C_1^{^{\frac{p}{{p'}}} } \left( {\frac{{2T_n }}{{2n + 1}}}
\right)^{\frac{p}{{p'}}} \frac{{2n + 1}}{{2T_n }}} \right)^{\frac{1}{p}}. $$
Since $T_n  = \frac{1}{2} + 3^{l - 1}  + n3^{l - 1} $,
we obtain  $\frac{{2T_n }}{{2n + 1}} \to
3^{l - 1} $ as $n \to \infty $, hence
$$\mathop {\overline {\lim } }\limits_{T \to \infty }
\left( {\frac{1}{{2T}}\int\limits_{ -
T}^T {|f(t) - g(t)|^p dt} } \right)^{\frac{1}{p}}
\ge C_2^{\frac{1}{p}}  \cdot 3^{\frac{1}{p}}
\cdot 3^{l\left( {\frac{1}{{p_0 }} - \frac{1}{p}} \right)}
- C_1^{\frac{1}{{p'}}}  \cdot
3^{\frac{1}{p} - \frac{1}{{p'}}}  \cdot 3^{l\left( {\frac{1}{{p'}}
- \frac{1}{p}} \right)}.$$
For  $p' > p_0 $, and $l$ sufficiently large the last inequality
contradicts to the equality
$$\mathop {\overline {\lim } }\limits_{T \to \infty } \left( {\frac{1}{{2T}}\int\limits_{ -
T}^T {\left| {f(t) - g(t)} \right|^p dt} } \right)^{\frac{1}{p}}  = d_{\left\{ 0 \right\}}^{B^p }
(f,g) = 0,$$
which is true for all $g(x)$ from the equivalent class of $f$
in the space $B_{\{ 0\} }^p $. The theorem is proved.

\end{document}